\documentclass{amsart}%
\usepackage{amsfonts}
\usepackage{amssymb}
\usepackage{amsmath}
\usepackage{graphicx}%
\providecommand{\U}[1]{\protect\rule{.1in}{.1in}}
\newtheorem{theorem}{Theorem}[section]
\newtheorem{lemma}[theorem]{Lemma}
\newtheorem{proposition}[theorem]{Proposition}

\newtheorem{example}[theorem]{Example}
\theoremstyle{remark}
\newtheorem{remark}[theorem]{Remark}
\numberwithin{equation}{section}

\newcommand{\abs}[1]{\lvert#1\rvert}
\newcommand{\R}{\mathbb{R}}
\newcommand{\C}{\mathbb{C}}
\newcommand{\Z}{\mathbb{Z}}
\newcommand{\D}{\mathbb{D}}

\newcommand{\T}{\mathbb{T}^1}

\newcommand{\ov}{\overline}
\newcommand{\begeq}{\begin{equation}}
\newcommand{\stopeq}{\end{equation}}
\newcommand{\ep}{\epsilon}
\newcommand{\dis}{\displaystyle}
\newcommand{\pa}{\partial}
\newcommand{\ei}[1]{\textrm{e}^{#1}}
\newcommand{\dd}[2]{\dis\frac{\pa #1}{\pa #2}}
\newcommand{\ta}{\theta}
\newcommand{\Si}{\Sigma}
\newcommand{\Ga}{\Gamma}
\newcommand{\Om}{\Omega}
\newcommand{\Omt}{\widetilde{\Omega}}
\newcommand{\til}{\widetilde}

\newcommand{\barom}{\ov{\Om}}
\newcommand{\siginf}{\Si^\infty}
\newcommand{\siga}{\Sigma^\infty_a}
\newcommand{\sigc}{\Sigma^\infty_c}
\newcommand{\Omsiga}{\Om\backslash\siga}
\newcommand{\Omjo}{\Om_j^0}
\newcommand{\kapj}{\kappa_j}
\newcommand{\cjk}{c_{jk}}
\newcommand{\Ajk}{A_{jk}}
\newcommand{\Gamjk}{\Ga_{jk}}
\newcommand{\pjkm}{p_{jk}^-}
\newcommand{\pjkp}{p_{jk}^+}
\newcommand{\Gamjone}{\Ga_{j1}}
\newcommand{\Gamjm}{\Ga_{jm}}
\newcommand{\Ajone}{A_{j1}}
\newcommand{\Ajm}{A_{jm}}
\newcommand{\varjk}{\vartheta_{jk}}
\newcommand{\aljk}{{\alpha_{jk}}}
\newcommand{\qjk}{q_{jk}}
\newcommand{\lam}{\lambda}
\newcommand{\lamj}{\lambda_j}
\newcommand{\Lam}{\Lambda}
\newcommand{\begar}{\begin{array}}
\newcommand{\stopar}{\end{array}}
\newcommand{\intc}[1]{\mathrm{Int}(C_{#1})}

\begin{document}

\author{A. Ainouz$^1$, K. E. Boutarene$^1$, A. Meziani$^2$}
\address{
\small $^1$Lab. AMNEDP, Department of Mathematics,
 USTHB,
 Po Box 32 El Alia,
Bab Ezzouar 16111
Algiers, Algeria}
\address{
\small $^2$Department of Mathematics, Florida International
University, 11200 S.W. 8th Street,  Miami, FL, 33199, USA e-mail: meziani@fiu.edu}

\title[ Riemann-Hilbert Problem ]{The Riemann-Hilbert problem for elliptic
vector fields with degeneracies}

\begin{abstract}
The boundary value problem
\[
Lu=f\ \ \text{in}\ \Om\, ,\quad \text{Re}(\ov{\Lam}\, u)=\phi\ \ \text{on}\ \pa\Om\, ,
\]
is studied for a class of planar complex vector fields $L$ in a simply connected open set $\Om\subset\R^2$.
The first integrals of $L$ are used to reduce the problem into a
collection of classical Riemann-Hilbert problems with discontinuous data.
\end{abstract}

\maketitle

\section*{Introduction}
This paper deals with boundary value problems for planar elliptic vector fields with degeneracies.
The vector field $L$ is nonsingular, complex valued, and real analytic in an open set containing
$\ov{\Om}$,
where $\Om$ is open, bounded, simply connected, and has a smooth boundary.
We assume that $L$ is locally solvable and is elliptic except on the analytic variety
$\Si$ given by $L\wedge\ov{L}=0$.
The boundary value problems considered here are versions of the Riemann-Hilbert problem
(RH problem for short) for the vector field $L$:
\begeq\left\{\begar{ll}
Lu =f & \quad\mathrm{in}\ \ \Om\, ,\\
\mathrm{Re}(\ov{\Lam}\, u)=\phi &\quad\mathrm{on}\ \ \pa\Om\, ,
\stopar\right.\stopeq
where $\Lam\, ,$ $\,\phi$ are H\"{o}lder continuous on $\pa\Om$ with $\abs{\Lam}=1$, $\phi$ is $\R$-valued,
and $f\in C^\infty(\ov{\Om})$.

The analytic variety $\Si$ consists of the set $\Si^0$, of points of finite type, and the set $\siginf$,
of points of infinite type (see Section 1 for definitions). $\siginf$ is a one-dimensional submanifold
in $\ov{\Om}$ and $L$ is tangent to $\siginf$. The connected components of $\siginf$ are the
(one-dimensional) orbits of $L$. We assume that these orbits are minimal (to be defined below).
Let $\Ga_1,\,\cdots ,\, \Ga_N$ be the non closed orbits  of $L$ such that each $\Ga_j$ connects in
$\Om$ two distinct points of $\pa\Om$. These orbits play a crucial role in the solvability
of Problem (0.1). For a given function $\Lam$, defined on $\pa\Om$,  we associate an index $\kapj\in\Z$ in each
connected component $\Om_j$ of $\Om\backslash\Ga_1\cup\cdots\cup\Ga_N$, and to each orbit
$\Ga_j$ with ends $p_j^-$, $p_j^+$ on $\pa\Om$, we attach a unique real number $\alpha_j\in [0,\ 1)$.
This number measures, in some sense, the jump of $\Lam$ along the orbit.

In general, it is not possible to find continuous solutions of Problem (0.1) on $\ov{\Om}$,
even in simple cases. The reason is that some of the indices $\kapj$ might be negative.
We will then allow for solutions to have isolated singular points: at most one singular point in each component
$\Om_j$ with index $\kapj <0$. In fact, the solution $u$ is continuous on
$\ov{\Om}\backslash\{p_1,\cdots ,p_n\}$ and, near each singular point $p_j$, it behaves
(through first integrals) as a meromorphic function with a pole of order $\le -\kapj$.
For the case $f=0$, the solutions $u$ are uniquely determined on the orbits $\Ga_j$ by
the values of $\Lam$ and $\phi$ at the ends $p_j^-$, $p_j^+$ and the number $\alpha_j >0$.
This value is given by the formula
\begeq
u(\Ga_j)=\ep\,\frac{\Lam(p_j^-)\phi(p_j^+)-\Lam(p_j^+)\phi(p_j^-)}{i\sin(\pi\alpha_j)}
\stopeq
with $\ep =\pm 1$.

Various properties of the types of vector fields considered here have been studied by many authors
and our approach is within the framework of the work contained in
{\cite{Ber-Mez2}},
{\cite{Ber-Cor-Hou}}, {\cite{Berh-Mez2}}, {\cite{Hou1}}, {\cite{Hou2}}, {\cite{Mez-JFA}},
{\cite{Mez-TAM}}, {\cite{Mez-MEM}}, {\cite{Nir-Tre}}, {\cite{Tre-Boo}}.
To our knowledge, boundary value problems for elliptic vector fields with interior degeneracies have not been studied,
except in {\cite{Mez-JDE}} where a special case of Problem (0.1) is considered when the
boundary functions satisfy $\Lam(p^-)=\Lam(p^+)$ and $\phi(p^-)=\phi(p^+)$ on each
orbit $\Ga$ with ends $p^-$ and $p^+$. This paper can thus be considered as an outgrowth of
Section 12 of {\cite{Mez-JDE}}.

Our basic technique to investigate Problem (0.1) is to make use of the first integrals of $L$
in the connected components of $\Om\backslash\Ga_1\cup\cdots\cup\Ga_N$ and to reduce the problem to
(N+1) classical RH problems for holomorphic functions in the unit disk $\D\subset\C$.
However, these new RH problems have discontinuities on $\pa\D$. These types of discontinuous RH problems
have been studied by many authors (see in particular the monographs {\cite{Beg}}, {\cite{Gak}}, and
{\cite{Mus}}). In this paper, we adapt the techniques contained in {\cite{Beg}}.

The organization of this paper is as follows. In Section 1, we recall the main properties of
the class of vector fields under consideration. In Section 2, we decompose the open set
$\Om\backslash\siginf$
into appropriate components and construct global first integrals for $L$. The notion of induced indices
for a boundary function $\Lam$ is introduced in Section 3. In Sections 4 and 5, we study Problem (0.1).
Theorem 4.1 gives the general continuous solution of Problem (0.1) when $f=0$ and $\phi =0$.
Theorem 5.2 deals with the case $f=0$ and $\phi\ne 0$ and gives the general (meromorphic) solution
satisfying (0.2) and in Theorem 5.4, the general case, with $f\ne 0$, is deduced.

\section{A class of vector fields}
We summarize some properties of the class of vector fields under consideration.
Most of the material listed in this section can be found in {\cite{Ber-Cor-Hou}},
{\cite{Mez-JDE}}, and {\cite{Tre-Boo}}.

Let $\Omt$ be an open subset of $\R^2$ and let $L$ be a vector field given by
\begeq
L=A(x,y)\dd{}{ x}+B(x,y)\dd{}{ y}\, ,
\stopeq
with $A,\, B\, \in C^\omega (\Omt,\C)$, i.e. $A$ and $B$ are real analytic and $\C$-valued in
$\Omt$. We assume that $L$ is nonsingular in $\Omt$ so that $\abs{A}+\abs{B} \, >0$ everywhere in
$\Omt$. We denote by $\ov{L}$ the complex conjugate vector field ( $\ov{L}=\ov{A}\pa_x+\ov{B}\pa_y$)
and view $L$ as a partial differential operator with symbol $\sigma =A\xi+B\eta$.
The operator $L$ is elliptic at a point $(x_0,y_0)$, if  $\sigma(x_0,y_0;\xi,\eta)\ne 0$ for every
$(\xi,\eta)\in\R^2\backslash 0$. The ellipticity of $L$ at $(x_0,y_0)$ is equivalent to $L$ and $\ov{L}$ being
independent at $(x_0,y_0)$. We denote by $\Sigma (\Omt)$ the set of points in $\Omt$ where $L$ fails to be elliptic.
 $\Sigma(\Omt)$ is the base projection of the characteristic set of $L$ (we often refer to $\Sigma(\Omt)$ as the characteristic
set of $L$). Consider the $\R$-valued function $\vartheta$ defined in $\Omt$ by
\begeq
\frac{1}{2i}L\wedge\ov{L} =\vartheta(x,y)\dd{}{x}\wedge\dd{}{y} =\mathrm{Im}(A(x,y)\ov{B(x,y)})\dd{}{x}\wedge\dd{}{y}\, .
\stopeq
Thus, $\Sigma(\Omt)$ is the analytic variety in $\Omt$ given by
\begeq
\Sigma(\Omt)=\left\{(x,y)\in\Omt\, ;\ \ \vartheta (x,y)=0\, \right\}\, .
\stopeq
We assume that $L$ is not a multiple of a real vector field, so that $\Sigma(\Omt)$ is a real analytic variety of dimension
$\le 1$.
A point $p_0\in\Sigma(\Omt)$ is of {\it{finite type}} if there exists a vector field $Y$ in the Lie algebra generated by
$L$ and $\ov{L}$ such that $L$ and $Y$ are independent at $p_0$. If no such vector field $Y$ exist, the point $p_0$
is said to be of {\it{inifinite type}}. The set of points $p\in\Si(\Omt)$ where $L$ is of finite type will be denoted
$\Si^0(\Omt)$ and those where $L$ is of infinite type will be denoted  $\Si^\infty(\Omt)$. We have then
\[
\Si(\Omt)=\Si^\infty(\Omt)\cup\Si^0(\Omt)\qquad\mathrm{and}\qquad
\Si^\infty(\Omt)\cap\Si^0(\Omt)=\emptyset\, .
\]
Since the vector field $L$ is nonsingular and real analytic, then it is {\it{locally integrable}}. This means
that for every point $p\in\Omt$, there exist an open set $\mathcal{O}_p\subset\Omt$, $\ p\in\mathcal{O}_p$, and
a function $z_p\in C^\omega(\mathcal{O}_p,\C)$ such that $Lz_p=0$ and $dz_p\ne 0$. The function $z_p$ is a local
first integral of $L$. The use of the local first integrals allows us to obtain local normal forms for
the vector field $L$.

In a neighborhood of an elliptic point $p\in\Omt\backslash\Si(\Omt)$, the vector field $L$ is conjugate
to the CR operator. That is, there exist local coordinates $(s,t)$ centered at $p$ in which $z_{p}=s+it$
is a first integral and $L$ is a multiple of $\dd{}{s}+i\dd{}{t}$. In a neighborhood of a point of finite type
$p_0\in\Si^0(\Omt)$, there exist local coordinates $(s,t)$ centered at $p$, an integer $k\in\Z^+$, an
$\R$-valued and $C^\omega$-function $\phi(s,t)$ satisfying
\begeq
\dd{^j\phi}{t^j}(0)=0,\ \ j=1,\cdots k-1,\ \ \dd{^k\phi}{t^k}(0)\ne 0
\stopeq
such that $z_p=s+i\phi(s,t)$ is a first integral and $L$ is conjugate to a multiple of
\[
\left(1+i\dd{\phi}{s}\right)\dd{}{t}-i\dd{\phi}{t}\dd{}{s}\, .
\]
Moreover, the variety $\Si^0(\Omt)$ is locally given by $\{\phi_t =0\}$.
In a neighborhood of a point of infinite type $p\in\Si^\infty(\Omt)$, there exist coordinates $(s,t)$ centered
at $p$, an integer $l\in\Z^+$, an $\R$-valued and $C^\omega$-function $\psi(s,t)$ with $\psi(0,t)\not\equiv 0$
such that $z_p=s+is^l\psi(s,t)$ is a first integral and $L$ is conjugate to a multiple of
\begeq
\left(1+i\dd{s^l\psi}{s}\right)\dd{}{t}-is^l\dd{\psi}{t}\dd{}{s}\, .
\stopeq
In these local coordinates, $\Si^\infty(\Omt)$ is the $t$-axis ($\{s=0\}$) and $\Si^0(\Omt)$ is given by
$\{\psi_t =0\}$ so that
\begeq
\Si(\Omt)=\{ (s,t);\ \ s^l\psi_t(s,t)=0\ \}\, .
\stopeq
It follows then that $\Si^\infty(\Omt)$ is a one-dimensional real analytic manifold embedded in $\Omt$
and that $L$ is tangent to $\Si^\infty(\Omt)$. For a given open set $U\subset\Omt$,
the connected components of $\Si^\infty(\Omt)\cap U$ will be refereed
to as the {\it{orbits}} of $L$ in $U$.

This paper deals with the vector fields with {\it minimal} orbits (as defined in {\cite{Mez-JDE}}).
That is, we will assume throughout that if $p\in\Si^\infty(\Omt)$, then the integer $l$ appearing in the
normalization (1.5) is $l=1$. The semiglobal normalization of vector fields with minimal orbits
is given in {\cite{Mez-JDE}}.

A vector field $L$ is said to be {\it locally solvable} at a point $p\in\Omt$, if there exists an
open set $\mathcal{O}\subset\Omt$, with $p\in\mathcal{O}$,  such that for every function $f\in C^\infty_0(\mathcal{O})$
the equation $Lu=f$ has a distribution solution $u\in\mathcal{D}'(\mathcal{O})$. A vector field $L$ is
solvable (and hypoelliptic) at all elliptic points $p\in\Omt\backslash\Si(\Omt)$. The solvability of $L$
at the nonelliptic points $p\in\Si(\Omt)$ is given by the Nirenberg-Treves Condition ($\mathcal{P}$).
In the case considered here of vector fields in two variables, this condition has a simple formulation.
The vector field $L$ is solvable at $p\in\Si$ if and only if there exists an open set $\mathcal{O}\subset\Omt$, with
 $p\in\mathcal{O}$,
such that the function $\vartheta =\mathrm{Im}(A\ov{B})$ defined in (1.2) does not change sign in any
connected component of $\mathcal{O}\backslash\siginf(\Omt)$. It follows, in particular, that if $L$ is locally
solvable at a point of finite type  $p\in\Si^0(\Omt)$, then the local first integral about $p$ is a local
homeomorphism (see {\cite{Tre-Boo}}). Furthermore, any continuous solution of $Lu=0$ in a
neighborhood of a point $p\notin\siginf(\Omt)$ can be written as $u=H\circ z$ where $z$ is a local first
integral and where $H$ is a holomorphic function defined in a neighborhood of $z(p)\in\C$. For solutions
with isolated singularities, we define a {\it{ pole of order}} $s$ of a solution $u$ to be a point
$p\notin\siginf$ such that there exists a meromorphic function $M$ with pole of order $s$ at $z(p)$ such that
$u=M\circ z$ in a neighborhood of $p$.
The global solvability of vector fields with minimal orbits is studied
in {\cite{Mez-JDE}}.

From now on, we will assume that $L$ satisfies Condition ($\mathcal{P}$) of local solvability and that it
has only minimal orbits.

\section{Decomposition of $\Om$ and global first integrals}

We decompose an open set $\Om$ into suitable components and construct global first integrals
for the vector field $L$. Let $L$, given by (1.1) with $A,B\, \in C^\omega(\Omt,\C)$, be
such that $L$ satisfies Condition ($\mathcal{P}$) and has only minimal orbits. Let $\Om\subset\R^2$
be a  simply connected and bounded open set with $\ov{\Om}\subset\Omt$. For simplicity, we will assume that
$\pa\Om$ is smooth (of class $C^\infty$) and that it is transversal to the manifold $\siginf(\Omt)$.
Throughout, we will use the following notation
\[
\Si =\Si(\Omt)\cap\ov{\Om}\, ,\quad \siginf=\siginf(\Omt)\cap\ov{\Om}\, \quad\mathrm{and}\quad
\Si^0=\Si^0(\Omt)\cap\barom\, .
\]
The set $\siginf$ is therefore a one-dimensional manifold. We denote by $\sigc$ the union of the closed curves
in $\siginf$ and by $\siga$ the union of the non-closed components of $\siginf$. We have then
\[
\siginf =\sigc\,\cup\,\siga\quad\mathrm{and}\quad\sigc\cap\siga =\emptyset\, .
\]
Hence, $\sigc$ consists of finitely many closed curves, each contained in the open set $\Om$, and $\siga$ consists of
finitely many arcs, each connecting two distinct points of $\pa\Om$ (since $\pa\Om$ is transversal to $\siginf(\Omt)$).

Let $\Ga_1\,,\cdots\, ,\ \Ga_N$ be the connected components of $\siga$. It is proved in {\cite{Mez-JDE}} that
\begeq
\Omsiga =\Om_1\cup\cdots\cup\Om_{N+1}\, ,
\stopeq
where each component is simply connected and for any two distinct components $\Om_j$ and $\Om_k$, we have
either $\barom_j \cap\barom_k =\emptyset$ or $\barom_j\cap\barom_k=\Ga_l$, for some component $\Ga_l$ of $\siga$.
The boundary $\pa\Om_j$ of each component $\Om_j$ of $\Omsiga$ is a piecewise smooth curve and consists of an arc
in $\pa\Om$ with ending points on $\siga$, followed by a component of $\siga$, followed by an arc in $\pa\Om$ and
so on. More precisely, if $\siga\ne\emptyset$, then for every $j=1,\cdots ,N+1$, there exists an integer $m\in\Z^+$
such that
\begeq
\pa\Om_j =\Ajone\cup\Gamjone\cup\cdots\cup\Ajm\cup\Gamjm\, ,
\stopeq
where the $\Ajk$'s are connected components of $\pa\Om\cap\pa\Om_j$ and the $\Gamjk$'s connected components
of $\pa\Om_j\cap\siga$. Moreover, $\Gamjk$ intersects transversally $\Ajk$ at a point $\pjkm$ and intersects
$A_{j(k+1)}$ at a point $\pjkp$. Hence,
\begeq
\Gamjk =\mathrm{arc}(\pjkm,\pjkp)\quad\mathrm{and}\quad \Ajk =\mathrm{arc}(p_{j(k-1)}^+,\pjkm)\, ,
\stopeq
with the understanding that $p_{j0}^+=p_{jm}^+$.
\begin{figure}[ht]
\flushleft 
\scalebox{0.45} {\includegraphics{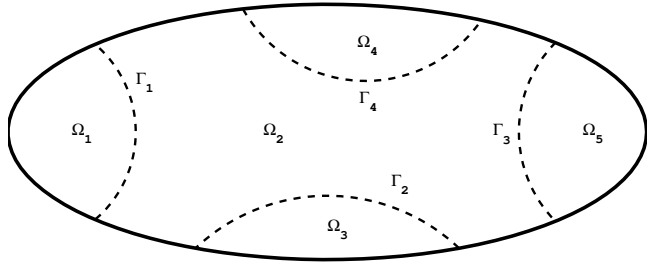}}
\caption{\textit{Decomposition of $\Omega$ by the $\Ga$-orbits}}
\end{figure}

For each $j\in\{1,\cdots ,N+1\}$, let
\[
\Om_j\backslash\sigc=\Omjo\cup\widehat{\Om}_j\,  ,
\]
where $\Omjo$ is the nonrelatively compact component of $\Om_j\backslash\sigc$ and $\widehat{\Om}_j$
is the union of the relatively compact components. Hence, there exist closed orbits $C_1,\,\cdots\, ,\,C_p$
of $L$ contained in $\Om_j$ such that
\begeq
\Om_j=\Omjo\cup\ov{\intc{1}}\cup\cdots\cup\ov{\intc{p}}\, ,
\stopeq
where $\intc{}$ denotes the open domain bounded by the closed curve $C$. It follows then that
\begeq
\pa\Om_j^0=\pa\Om_j \cup C_1\cup\cdots\cup C_p\, .
\stopeq
The following result is proved in {\cite{Mez-JDE}}.

\begin{figure}[ht]
\flushleft
\scalebox{0.55} {\includegraphics{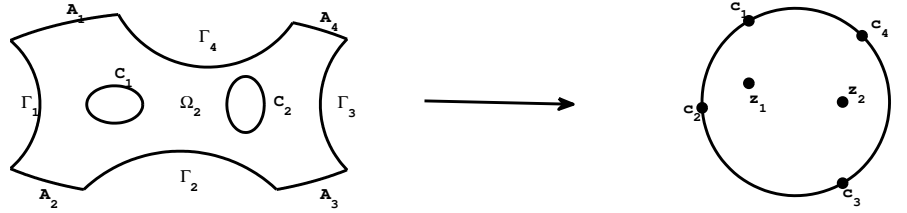}}
\caption{\textit{First integral $Z$ on $\Omega_2$}}
\end{figure}

{\proposition{\rm ({\cite{Mez-JDE}})} {Let $\Om_j$ and $\Omjo$ be as above with boundaries
given by {\rm (2.2)} and {\rm (2.4)}, respectively. Then there exist a function
\[
F_j:\ \ov{\Omjo}\, \longrightarrow\ \, \ov{\mathbb{D}}\, ,
\]
where $\mathbb{D}$ is the unit disk in $\C$; points $z_{j1},\cdots ,z_{jp}\in\mathbb{D}$;
and points $c_{j1},\cdots ,c_{jm}\in\pa\mathbb{D}$ such that:
\begin{itemize}
\item $F_j$ is $C^\omega$ in $\Omjo$ and H\"{o}lder continuous on $\ov{\Omjo}$;
\item $LF_j=0$ in $\Omjo$ and $dF_j(p)\ne 0$ for every $p\in\Omjo$;
\item $F_j:\ \Omjo\,\longrightarrow\, \mathbb{D}\backslash\{z_{j1},\cdots,z_{jp}\}$ is a homeomorphism;
\item $F_j(C_l)=z_{jl} $ for $l=1,\cdots ,p$ and $F_j(\Gamjk)=c_{jk}$ for $k=1,\cdots ,m$.
\end{itemize}
}}

{\remark{ Choose an orientation of $\R^2$ and its induced orientation on the components $\Om_j$'s of
$\Omsiga$. Given a component $\Om_j$, if the first integral $F_j$ preserves the orientation, then any other
first integral also preserves the orientation. The components $\Om_j$'s can thus be divided into two
types: Those with orientation preserved  and those with orientation reversed by first integrals.
Furthermore, if $\Om_j$ and $\Om_k$ are two adjacent components ($\ov{\Om_j}\,\cap\,\ov{\Om_k}\ne\emptyset$),
then one of the component is orientation preserved and the other orientation reversed (see {\cite{Mez-JDE}}
for details).  } }

\medskip

The first integrals $F_j$ of Proposition 2.1 can be glued together to obtain a global first integral in
\begeq
\Om^0=\Om_1^0\cup\cdots\cup\Om_{N+1}^0\, .
\stopeq
More precisely, we have the following theorem.

{\theorem{There exists a function $Z:\, \ov{\Om^0}\,\longrightarrow\,\ov{\mathbb{D}}$,
where $\Om^0$ is given by {\rm (2.6)}, such that:
\begin{itemize}
\item $Z\in C^\omega(\Om^0,\mathbb{D})$, $\ Z$ is H\"{o}lder continuous on $\ov{\Om^0}$;
\item $LZ=0$ and $dZ\ne 0$ in $\Om^0$;
\item For each $j=1,\cdots ,N+1$, there exist points $z_{j1},\cdots , z_{jp}\in\mathbb{D}$
such that
\[
Z(C_l)=z_{jl}\ \mathrm{for}\ l=1,\cdots , p\ \ \mathrm{and}\ \
Z: \ \Om_j^0\, \longrightarrow\, \mathbb{D}\backslash\{z_{j1},\cdots ,z_{jp}\}
\]
is a homeomorphism;
\item For each $j=1,\cdots ,N+1$ and $k=1,\cdots ,m$ there exist $\cjk\in\pa\mathbb{D}$
such that
\[
Z(\Gamjk)=\cjk\quad\mathrm{and}\quad Z(\Ajk)=\mathrm{arc}(c_{j(k-1)},\cjk)\, .
\]
\end{itemize}}}

\begin{proof}
Select a component $\Om_j$ of $\Omsiga$ and label it $V_0$ so that, as in (2.4), we have
\[
V_0=V_0^0\cup\ov{\intc{1}}\cup\cdots\cup\ov{\intc{p}}\, ,
\]
where $V_0^0$ is the non relatively compact component of $V_0\backslash\siga$, and with boundary
\[
\pa V_0=A_1\cup\Ga_1\cup\cdots\cup A_m\cup\Ga_m\, ,
\]
with the $A_l$'s and $\Ga_l$'s as in (2.2) and (2.3). Let $F^0:\, \ov{V_0^0}\longrightarrow\ov{\mathbb{D}}$
be a first integral of $L$ in $V_0^0$ satisfying Proposition 2.1. Set, as above,
$\Ga_k=\mathrm{arc}(p_k^-,p_k^+)$, $\ A_k=\mathrm{arc}(p_{k-1}^+,p_k^-)$, $\ F^0(C_s)=z_s\in\mathbb{D}$,
and $F^0(\Ga_k)=c_k\in\pa\mathbb{D}$ so that $F^0(A_k)=\mathrm{arc}(c_{k-1},c_k)$.

For each $k=1,\cdots , m$ there exists a unique component $U_k$ of $\Omsiga$ such that
$\ov{V^0}\cap\ov{U_k}=\Ga_k$. Let $F_k:\ov{U_k^0}\longrightarrow\ov{\mathbb{D}}$ be a first integral of
$L$ in $U_k^0$ as in Proposition 2.1. If $F_k(\Ga_k)=d_k\in\pa\mathbb{D}$, let
$\til{F_k}:\ov{U_k^0}\longrightarrow\ov{\mathbb{D}}$ be given by
\[
\til{F_k}(x,y)=\frac{c_k}{d_k}F_k(x,y).
\]
$\til{F_k}$ is another first integral of $L$ in $U_k^0$ satisfying Proposition 2.1 and
$\til{F_k}(\Ga_k)=F^0(\Ga_k)=c_k$. Hence $\til{F_k}$ is a H\"{o}lder continuous extension of
$F^0$ to $U_k^0$ (see {\cite{Mez-JDE}}). Let
\[
V_1=V_0\cup U_1\cup\cdots \cup U_m\, .
\]
Define
\[
F^1:\ \ov{V_1^0}=\ov{V_0^0}\cup\ov{U_1^0}\cup\cdots\cup\ov{U_m^0}\,\longrightarrow\,
\ov{\mathbb{D}}
\]
by
\[
F^1(x,y)=\left\{\begar{ll}
F^0(x,y) & \mathrm{if}\quad (x,y)\in \ov{V_0^0};\\
\til{F_k}(x,y) & \mathrm{if}\quad (x,y)\in \ov{U_k^0}\, .
\stopar\right.
\]
Then, $F^1$ is a first integral of $L$ in $V_1^0$ that satisfies the properties of the theorem.
This extension of the first integral can be repeated through each orbit contained in $\pa V_1$ to a larger union of
components $V_2$ and so on until we exhaust all the components of $\siga$ and reach the global first
integral in $\Om^0$.
\end{proof}

\section{Induced indices}
For a function defined on the boundary of $\Om$, we define its indices with respect to the
components of $\Omsiga$. Let $L$ and $\Om$ be, respectively, a vector field and an open set as in Section 2.
Thus, $\Om$ is simply
connected, $\ov{\Om}\subset\Omt$, $\ \pa\Om$ is a smooth closed curve and transversal to $\siga$.
We assume that $\pa\Om$ is positively oriented (counterclockwise). Let
\begeq
\Omsiga =\Om_1\cup\cdots\cup\Om_{N+1}
\stopeq
be the decomposition given in Section 2. We assume that $N\ge 1$ and that $\Om_1$ is
orientation preserved by the first integrals of $L$. Fix a point
$s_0\in\pa\Om_1\backslash\siga$ and let $\tau : [0,\ 2\pi ]\longrightarrow\pa\Om$ be
a smooth parametrization such that $\tau(0)=\tau(2\pi)=s_0$.

Let $\Lam :\,\pa\Om\,\longrightarrow\,\T =\pa\mathbb{D}$ be a H\"{o}lder continuous
function ($\Lam\in C^\sigma (\pa\Om,\T)$ with $0<\sigma<1$) and denote by $\arg\Lam$ a
continuous branch of the argument of $\Lam$. That is, $\Lam(s)=\ei{i\arg\Lam(s)}$
and $\arg\Lam$ is defined through the parametrization $\tau$ as a continuous function
\[
\arg\Lam\circ\tau :\, [0,\ 2\pi]\,\longrightarrow\, \R
\]
such that
\begeq
\arg\Lam(s_0^-)=\arg\Lam(s_0^+)+\pi q_0\quad\mathrm{with}\quad q_0\in\Z\, ,
\stopeq
where $\dis\arg\Lam(s_0^+)=\lim_{t\to 0^+}\arg\Lam(\tau(t))$ and
$\dis\arg\Lam(s_0^-)=\lim_{t\to 2\pi^-}\arg\Lam(\tau(t))$. Note that
$q_0/2$ is the winding number of $\Lam$.

Now  we define the jump along an orbit $\Ga\subset \siga$ of the function $\Lam$.
Such an orbit $\Ga$ is the common boundary of two adjacent components $\Om_{j}$ and
$\Om_{l}$ of $\Omsiga$. One of the components, say $\Om_{j}$ is orientation preserved by
the first integrals and the other is orientation reversed (Remark 2.2). Consider then $\Om_{j}$
and its boundary
\begeq
\pa\Om_j =\Ajone\cup\Gamjone\cup\cdots\cup\Ajm\cup\Gamjm
\stopeq
as in (2.2). The orbit appears as $\Ga =\Gamjk$ for some $k$. The boundary $\pa\Om_j$ inherits its orientation
from that of $\pa\Om$ and so do  each of its arcs
$A_{jl}$ and $\Ga_{jl}$.
Let $\pjkm$ and $\pjkp$ be the ordered ends of $\Gamjk$. We have
\begeq
\Gamjk=\mathrm{arc}(\pjkm,\pjkp)\quad\mathrm{and}\quad
\Ajk=\mathrm{arc}(p_{j(k-1)}^+,\pjkm)\, .
\stopeq
We define the jump of $\Lam$ along $\Ga=\Gamjk$ as the real number
\begeq
\varjk=\arg\Lam(\pjkm)-\arg\Lam(\pjkp)\, .
\stopeq

To define the index of $\Lam$ relative to the component $\Om_j$, we start by writing the jumps
$\varjk$ as
\begeq
\varjk=\pi\aljk +\pi\qjk\, ,
\stopeq
with
\begeq
\qjk=\left[\frac{\varjk}{\pi}\right]\, \in\Z\quad\mathrm{and}\quad
\aljk=\frac{\varjk}{\pi}-\qjk \,\in [0,\ 1)\, ,
\stopeq
where, for a real number $r$, $\ [r]$ denotes the largest integer $\le r$.
Now, we distinguish two types of orbits $\Gamjk$. Those for which the integer $\qjk$
is even and those for which it is odd. Let
\begeq
\mathcal{C}_j^1=\left\{\Gamjk;\ \qjk\in 2\Z+1\right\}\, ,\qquad
\mathcal{C}_j^2=\left\{\Gamjk;\ \qjk\in 2\Z\right\}\, ,
\stopeq
and let $n_j$ be the number of elements in $\mathcal{C}_j^1$. We define the index $\kapj$ of
$\Lam $ relative to $\Om_j$ as follows:
\begeq
\begar{ll}
\dis\kapj=\mathrm{Ind}(\Lam,\Om_j) =\frac{1}{2}\left(\sum_{k=1}^m\qjk\, -n_j\right) &\ \mathrm{if}\
j=2,\cdots ,N+1\, , \\
\dis\kappa_1=\mathrm{Ind}(\Lam,\Om_1) =\frac{1}{2}\left(q_0+\sum_{k=1}^mq_{1,k}\, -n_1\right)\, . &
\stopar
\stopeq
Note that $\kapj\in\Z$ for $j=1,\cdots ,N+1$.

{\example\rm{ Let
\begeq
L=(1+2ixy)\dd{}{y}-i(x^2-1)\dd{}{x}\, .
\stopeq
It can be seen at once that $L\wedge\ov{L}=-2i(x^2-1)\dis\dd{}{x}\wedge\dd{}{y}$. Hence, $L$ is elliptic
everywhere in $\R^2$ except on the vertical lines $x=1$ and $x=-1$ along which $L$ is tangent.
We have then $\Si=\siga=\{x=\pm 1\}$. Furthermore, the orbits $x=1$ and $x=-1$ are minimal.
Note that the function $F(x,y)=x+iy(x^2-1)$ is a global first integral of $L$ in $\R^2$.

Let $\Om=D(0,2)$ be the disk with center $0$ and radius $2$ in $\R^2$. Then
\[\Omsiga=\Om_1\cup\Om_2\cup\Om_3\, ,\]
where
\begeq\begar{ll}
\Om_1 & =\left\{ (x,y)\in\R^2:\ \ 1<x,\ \  x^2+y^2<4\right\}\, ;\\
\Om_2 & =\left\{ (x,y)\in\R^2:\ \ -1<x<1,\ \  x^2+y^2<4\right\}\, ;\\
\Om_3 & =\left\{ (x,y)\in\R^2:\ \ x<-1,\ \  x^2+y^2<4\right\}\, .
\stopar\stopeq
$\siga$ consists of the  vertical segments
\[
\Ga_1=\left\{ (1,y):\ -\sqrt{3}\le y\le\sqrt{3}\right\}\quad\mathrm{and}\quad
\Ga_2=\left\{ (-1,y):\ -\sqrt{3}\le y\le\sqrt{3}\right\}\, .
\]
Since the map $\Phi=(\mathrm{Re}(F),\mathrm{Im}(F))=(x,y(x^2-1))$ is orientation preserving for $\abs{x}>1$
and orientation reversing for $\abs{x}<1$, then $\Om_1$ and $\Om_3$ are orientation preserved
by first integrals and $\Om_2$ is orientation reversed.
As oriented arcs we have $\Ga_1=[2\ei{i\pi/3},\ 2\ei{5i\pi /3}]$ and $\Ga_2=[2\ei{4i\pi/3},\ 2\ei{2i\pi /3}]$.

Let $\Lam(2\ei{i\ta})=\ei{if(\ta)}$ with $f:\, [0,\ 2\pi ]\longrightarrow\R$ smooth and
$f(2\pi)=f(0)+\pi q_0$ with $q_0\in 2\Z$. The jumps of $\Lam$ along $\Ga_1$ and $\Ga_2$ are
respectively
\[
\vartheta_1=f(\pi /3)-f(5\pi/3)\quad\mathrm{and}\quad
\vartheta_2=f(4\pi /3)-f(2\pi/3)\, .
\]
The induced indices are respectively
\[
\kappa_1=\frac{1}{2}(q_0+q_1-n_1),\quad\kappa_2=\frac{1}{2}(q_1+q_2-n_2)\quad\mathrm{and}
\quad\kappa_3=\frac{1}{2}(q_2-n_3)\, ,
\]
where $q_1=[\vartheta_1/\pi]$; $\ q_2=[\vartheta_2/\pi]$; $\ n_1=1$ if $q_1$ is odd and $0$ if not;
$\ n_2=2$ if both $q_1$ and $q_2$ are odd, $\ n_2=1$ if only one $q_1$ or $q_2$ is odd  and $n_2=0$ if both
are even; and $n_3=1$ if $q_2$ is odd and $0$ if not.

In the case $f(\ta)=\ta$, we have $q_0=2$, $\ \vartheta_1=-4\pi /3$, $\ \vartheta_2=2\pi/3$,
$\ q_1=-2$, $\ q_2=0$, $\ n_1=0$, $\ n_2=0$ and $\ n_3=0$. The  numbers $\aljk$ are in this case
$\alpha_1=2/3$ and $\alpha_2 =2/3$. The indices are $\kappa_1=0$, $\kappa_2=-1$, and $\kappa_3=0$

In the case $f(\ta)=\pi\sin\ta$, we have $q_0=0$, $\ \vartheta_1=\pi\sqrt{3}$, $\ \vartheta_2=-\pi\sqrt{3}$,
$\ q_1=1$, $\ q_2=-2$, $\ n_1=1$, $\ n_2=1$, $\ n_3=0$, $\ \alpha_1=\sqrt{3}-1$, $\ \alpha_2=2-\sqrt{3}$. The indices
are $\kappa_1=0$, $\ \kappa_2=-1$, and $\kappa_3=-1$.
}}

\section{The Homogeneous Riemann-Hilbert problem for $L$}
We consider the simplest boundary value problem associated with a vector field.
As in the previous section, $L$ is a nonsingular complex vector field with real analytic coefficients given by
(1.1) in an open set $\Omt\subset\R^2$. We assume that $L$ satisfies Condition ($\mathcal{P}$) and that it
has only minimal orbits. Let $\Om$ be a bounded, simply connected open set, with $\ov{\Om}\subset\Omt$ and with $ \pa\Om$
smooth and transversal to $\siginf(\Omt)$. As in the previous sections, we use the decompositions
\begeq
\Omsiga =\Om_1\cup\cdots\cup\Om_{N+1}
\stopeq
and for $j=1,\cdots ,N+1$,
\begeq
\pa\Om_j=\Ajone\cup\Gamjone\cup\cdots\cup\Ajm\cup\Gamjm
\stopeq
 where $\Gamjk$ and $\Ajk$ are as described  in Sections 2 and 3.
Let $\Lam\in C^\sigma(\pa\Om ,\pa\D)$ with $0<\sigma <1$. For each $j\in\{1,\cdots ,N+1\}$,
let $\kapj$ be the index of $\Lam$ relative to $\Om_j$ and for $k\in\{1,\cdots ,m\}$, let
$\aljk\in [0,\ 1)$ be the associated jump with the orbit $\Gamjk$ as defined in (3.7).
Denote by $\delta_j$ the number of orbits $\Gamjk$ with $\aljk =0$.
Consider the boundary value problem
\begeq\left\{\begar{ll}
Lu=0 &\quad\mathrm{in}\ \ \Om\, ,\\
\mathrm{Re}(\ov{\Lam} u)=0 &\quad\mathrm{on}\ \ \pa\Om\, .\stopar\right.\stopeq
We have the following theorem.

{\theorem{ Problem {\rm{(4.3)}} has nontrivial continuous solutions on $\ov{\Om}$ that
vanish on $\siga$ if and only if there exists $j\in\{1,\cdots ,N+1\}$ with
$2\kapj\ge\delta_j$. Furthermore, the number of independent such solutions is
\begeq
r=\sum_{ \kapj\ge\delta_j}\left(2\kapj-\delta_j+1\right)\, .
\stopeq}}

\begin{proof} For $j=1,\cdots ,N+1$, let $Z_j:\, \ov{\Om_j^0}\,\longrightarrow\,\ov{\D}$ be
the restriction to $\Om_j^0$ of a first integral $Z$ of $L$ satisfying Theorem 2.3. As in Section
2, $\Om_j^0$ is the  component of $\Om_j\backslash\sigc$ given in (2.4) and
$\pa\Om_j^0$ given by (2.5). The boundary $\pa\Om_j$ is as in (4.2) and
$\Gamjk=\mathrm{arc}(\pjkm,\pjkp)$, $\ \Ajk=\mathrm{arc}(p_{j(k-1)}^+,\pjkm)$
(see Section 3), set
\begeq
Z_j(\Gamjk)=\cjk\, \in\pa\D\quad\mathrm{and}\quad Z_j(\Ajk)=\mathrm{arc}(c_{j(k-1)},\cjk)
\,\subset\pa\D
\stopeq
(with $c_{j0}=c_{jm}$). For each $l=1,\cdots , p$, let $z_l=Z_j(C_l)\in\D$, where $C_l$
is a closed orbit of $L$ contained in $\pa\Om_j^0$ as in (2.5).

Define the function $\lamj$ on $\pa\D\backslash\{c_{j1},\cdots ,c_{jm}\}$
by $\lamj(\zeta)=\Lam\circ Z_j^{-1}(\zeta)$. Note that $\lamj$ is H\"{o}lder continuous on each arc
$\ov{Z_j(\Ajk)}\subset\pa\D$. The function $w_j=u\circ Z_j^{-1}$ satisfies the RH-problem
\begeq
\dd{w_j}{\ov{z}}=0\quad\mathrm{in}\ \D\, ,\qquad\mathrm{Re}(\ov{\lamj}\, w_j)=0\quad\mathrm{on}\ \pa\D\, .
\stopeq
This is a discontinuous RH-Problem that was considered by several authors (see for instance {\cite{Beg}},
{\cite{Gak}}, and {\cite{Mus}}). Note that with our definition of the jump $\varjk$ of $\Lam$ along $\Gamjk$,
we have
\begeq
\lamj(\cjk^-)=\ei{i\varjk}\lamj(\cjk^+)\, ,
\stopeq
where $\dis\lamj(c^\pm)=\lim_{t\to 0^\pm}\lamj(c\ei{it})$. Note also that (4.7) is valid whether
$Z_j$ is orientation preserving or not.

Now we proceed to the construction of solutions of Problem (4.6) as it is done in Chapter 3 of {\cite{Beg}}.
We transform the problem into a continuous RH-Problem and use the Schwarz operator to write the general solution.
We start by defining the function $\til{\lamj}$ on $\pa\D$ by
\begeq
\til{\lamj}(\zeta)=\lamj(\zeta)\prod_{k=1}^m\frac{\ov{(\zeta-\cjk)}^{\aljk}}{\abs{\zeta-\cjk}^{\aljk}}
\stopeq
($z^\aljk$ is defined through any branch of the logarithm). It is verified at once (see {\cite{Beg}}) that
\begeq
\til{\lamj}(\cjk^-)=(-1)^{\qjk}\til{\lamj}(\cjk^+)\, ,
\stopeq
where $\qjk$ is the integer associated with the orbit $\Gamjk$ and defined through the jump $\varjk$ by (3.7).
Hence, $\lamj$ is continuous at the points $\cjk$ when $\qjk$ is even. To obtain a continuous function through
the points with $\qjk$ odd, we consider the collection $\mathcal{C}_j^1=\{ \cjk:\ \qjk\ \mathrm{odd}\}$ (if not empty)
as a set of ordered points in $\pa\D$ (by a the argument function) and write it as
\begeq
\mathcal{C}_j^1=\{ c_{j\mu_s}:\ \ 1\le s\le a\}\, ,\quad \mu_s <\mu_{s+1} \quad (1\le s <a\le m)\, .
\stopeq
We distinguish two cases depending on whether the number of elements, $a$, in $\mathcal{C}_j^1$ is even or odd.
In case $a$ is even, define the arcs in $\pa\D$ with ends on $\mathcal{C}_j^1$ as
\begeq\begar{l}
\dis B_s =\{\zeta\in\pa\D :\ \ \arg c_{j\mu_s}<\arg\zeta<\arg c_{j\mu_{s+1}}\}\, ,\quad 1\le s <a\, ,\\

\dis B_a =\{\zeta\in\pa\D :\ \ \arg c_{j\mu_a}<\arg\zeta<\arg c_{j\mu_1}+2\pi\}\, .
\stopar\stopeq
Let $\beta:\, \pa\D\backslash\mathcal{C}_j^1\,\longrightarrow\,\{-1, 1\}$ be the alternating function
defined by
\begeq
\beta(\zeta)=(-1)^{s-1}\quad\mathrm{if}\quad \zeta\in B_s\quad 1\le s\le a\, .
\stopeq
In the case when $a$ is odd, consider an additional point $c_0\in\pa\D\backslash\{c_{j1},\cdots ,c_{jm}\}$
and replace $\mathcal{C}_j^1$ by $\widehat{\mathcal{C}_j^1}=\mathcal{C}_j^1\cup\{c_0\}$ and proceed to define
the arcs $B_s$'s and the alternating function $\beta$ accordingly.

Now we can define the continuous function $\lamj^0$ on $\pa\D$ by
\begeq
\lamj^0(\zeta)=\beta(\zeta)\lamj(\zeta)\prod_{k=1}^m\frac{\ov{(\zeta-\cjk)}^{\aljk}}{\abs{\zeta-\cjk}^{\aljk}}\quad
\mathrm{if}\ a\in 2\Z\, ,
\stopeq
\begeq
\lamj^0(\zeta)=\beta(\zeta)\lamj(\zeta)\frac{\ov{\zeta -c_0}}{\abs{\zeta-c_0}}\prod_{k=1}^m\frac{\ov{(\zeta-\cjk)}^{\aljk}}{\abs{\zeta-\cjk}^{\aljk}}\quad
\mathrm{if}\ a\in 2\Z +1\, .
\stopeq
It can be verified ({\cite{Beg}}) that such a function $\lamj^0$ is H\"{o}lder continuous on $\pa\D$
and that it has index $\kapj$.

Assume that $a$ is even, so that the function $\lamj^0$ is given by (4.13) (the case $a$ odd can be dealt with
in a similar fashion). Choose the branch of the logarithm so that $\arg(\zeta-\cjk)\,\to\, 0$ as
$\zeta\,\to \,\cjk^+$ on $\pa\D$ and $\arg(\zeta-\cjk)\,\to\,\pi$ as $\zeta\,\to\,\cjk^-$ on $\pa\D$. Then we
have for $k=1,\cdots ,m$
\begeq
\lamj^0(\cjk)=\beta(\cjk^+)\lamj(\cjk^+)\prod_{l\ne k}
\frac{\ov{(\cjk-c_{jl})}^{\alpha_{jl}}}{\abs{\cjk-c_{jl}}^{\alpha_{jl}}}\, .
\stopeq
Consider the function $w_j^0$ defined in $\D$ by
\begeq
w_j^0(z)=\frac{w_j(z)}{\prod_{k=1}^m(z-\cjk)^\aljk}\, .
\stopeq
The RH-Problem for $w_j^0$ is then
\begeq
\dd{w_j^0}{\ov{z}}=0\quad\mathrm{in}\ \D\, ,\qquad \mathrm{Re}(\ov{\lamj^0}\, w_j^0)=0\quad\mathrm{on}\ \pa\D\, .
\stopeq
This is the classical RH-Problem with only the trivial solution, if $\kapj <0$ and with general solution for
$\kapj\ge 0$, given by
\begeq
w_j^0(z)=z^{\kapj}\ei{i\gamma(z)}\left(
id_0+\sum_{l=1}^{\kapj}d_lz^l-\ov{d_l}\frac{1}{z^l}\right)\, ,
\stopeq
where $d_0\in\R$, $\ d_1,\cdots ,d_{\kapj}\in\C$ are arbitrary constants and where
\begeq
\gamma(z)=\mathcal{S}(\arg(\zeta^{-\kapj}\lamj^0(\zeta)))(z)\, ,
\stopeq
with $\mathcal{S}$ the Schwarz operator on the disk $\D$:
\begeq
\mathcal{S}(f)(z)=\frac{1}{2\pi i}\int_{\pa\D}f(\zeta)\,\frac{\zeta+z}{\zeta-z}\,\frac{d\zeta}{\zeta}\, .
\stopeq
Note that for $z=\zeta\in\pa\D$, we have $\mathrm{Re}(\gamma)(\zeta)=\arg\lamj^0(\zeta)-\kapj\arg\zeta$.

It follows from (4.16) and (4.18) that when $\kapj\ge 0$, the general solution
\begeq
w_j(z)=w_j^0(z)\prod_{k=1}^m(z-\cjk)^\aljk
\stopeq
of the Problem (4.6), vanishes at each point $\cjk$ with $\aljk >0$. At a point $c_{jk_0}\in\pa\D$
with $\alpha_{jk_0}=0$, the function $w_j$ vanishes only when the coefficients $d_0$, $d_1,\cdots, d_{\kapj}$
satisfy the linear equation
\begeq
id_0+\sum_{l=1}^{\kapj}d_lc_{jk_0}^l-\ov{d_lc_{jk_0}^l} =0\, .
\stopeq
Since there are $\delta_j$ points $c_{jk_0}$ (with $\alpha_{jk_0}=0$), then the ($2\kapj +1$)-real coefficients
$d_0$, $\ \mathrm{Re}(d_1)$, $\mathrm{Im}(d_1)$, $\,\cdots\,$,   $\mathrm{Re}(d_{\kapj})$, $\mathrm{Im}(d_{\kapj})$
must satisfy $\delta_j$ linear equations (over $\R$). We get then $2\kapj-\delta_j+1$ independent solutions
if $2\kapj\ge\delta_j$ and only the trivial solution, if $2\kapj<\delta_j$.

Finally, we can define the continuous solutions of the original Problem (4.3) in $\Om$ by
\begeq
u(x,y)=\left\{\begar{ll}
0 & \ \mathrm{if}\quad (x,y)\in\ov{\Om_j}\ \ \mathrm{with}\ \ 2\kapj <\delta_j\, ;\\
w_j(Z_j(x,y)) & \ \mathrm{if}\quad (x,y)\in\ov{\Om_j^0}\ \ \mathrm{with}\ \ 2\kapj \ge \delta_j\, ;\\
w_j(z_l) & \ \mathrm{if}\quad (x,y)\in\ov{\intc{l}}\subset{\Om_j}\ \ \mathrm{with}\ \ 2\kapj \ge \delta_j\, ,
\stopar\right.\stopeq
 where $w_j$ is given by (4.21) and the coefficients satisfy (4.22). It follows also that the number of
 independent solutions
 is given by (4.4).
\end{proof}

{\example \rm{We continue with the vector $L$ and $\Om=D(0,2)$ of Example 3.1. We have seen that
$\siga=\Ga_1\cup\Ga_2$ and $\Omsiga=\Om_1\cup\Om_2\cup\Om_3$.

For $\Lam(2\ei{i\ta})=\ei{i\ta}$, we found that $\kappa_1=\kappa_3=0$, $\, \kappa_2=-1$, and
$\delta_1=\delta_2=\delta_3=0$. Hence the RH-Problem
\[
Lu=0\quad\mathrm{in}\ D(0,2),\qquad\mathrm{Re}\left(\ei{-i\ta}u(2\ei{i\ta})\right)=0\quad
\mathrm{on}\ \pa D(0,2)
\]
has 2 independent solutions that are identically zero in $\Om_2$.

For $\Lam(2\ei{i\ta})=\ei{i\sin\ta}$, this time we have $\kappa_1=0$, $\,\kappa_2= \kappa_3=-1$, and
$\delta_1=\delta_2=\delta_3=0$. Hence the RH-Problem
\[
Lu=0\quad\mathrm{in}\ D(0,2),\qquad\mathrm{Re}\left(\ei{-i\sin\ta}u(2\ei{i\ta})\right)=0\quad
\mathrm{on}\ \pa D(0,2)
\]
has 1 independent solution that is identically zero in $\Om_2\cup\Om_3$.}}

\section{Nonhomogeneous RH problems}
In this section we consider nonhomogeneous RH-Problems for the vector field $L$ on a
simply connected open set $\Om$. The assumptions on $L$ and $\Om$ are as in Section 4.
The first problem to be studied is
\begeq\left\{\begar{ll}
Lu=0 & \ \ \mathrm{in}\ \ \Om\, ,\\
\mathrm{Re}\left(\ov{\Lam}\, u\right) =\phi &\ \ \mathrm{on}\ \ \pa\Om\, ,
\stopar\right.\stopeq
where the functions $\Lam$ and $\phi$ are H\"{o}lder continuous:
$\Lam\in C^\sigma(\pa\Om,\C)$ and $\phi\in C^\sigma(\pa\Om,\R)$, for some
$0<\sigma<1$. In general such a problem does not have continuous solutions
throughout $\Om$, even in most simple cases, due to the fact that some of the induced
indices of $\Lam$ might be negative. We will allow then for solutions to have a
single pole in each of the regions
with negative index (see Section 1 for the definition of a pole).
We will also  consider only the {\it generic} case when $\Lam$ takes different values
at the ends of each arc $\Ga\in\siga$ (so that the corresponding number $\alpha$
defined in (3.7) is positive: $0<\alpha <1$).  Before, we state the main result about
Problem 5.1, we recall the following property of Cauchy type integrals that will be used
and  whose proof can be found in {\cite{Mus}} page 85.

{\lemma {{\rm ({\cite{Mus}})}
Given a positively oriented curve $\mathcal{A}=\mathrm{ arc}(A,B)$ in $\C$; a point
$c\in\mathcal{A}$; a function $g$ on $\mathcal{A}$ that  is H\"{o}lder continuous
on  $\ov{\mathrm{arc}(A,c)}$ and on $\ov{\mathrm{arc}(c, B)}$; and given a real number $0<\alpha <1$,
then the Cauchy integral
\begeq
\Phi(z)=\frac{1}{2\pi i}\int_{\mathcal{A}}\,\frac{g(\zeta)}{(\zeta-z)^\alpha (\zeta -z)}\, d\zeta
\stopeq
satisfies the following:
\begeq
\Phi(z)=\frac{\ei{i\pi\alpha}g(c^+)-\ei{-i\pi\alpha}g(c^-)}{2i\sin(\pi\alpha)}\,\frac{1}{(z-c)^\alpha}+
\Phi_0(z)\, ,
\stopeq
for $z\in\C\backslash\mathcal{A}$ to the left of $\mathcal{A}$, with $\Phi_0$ analytic and such that
\begeq
\abs{\Phi_0(z)}\le\,\frac{K}{\abs{z-c}^{\alpha_0}}
\stopeq
for some $K>0$ and $0<\alpha_0<\alpha$.
}}

\medskip

We have the following theorem.

{\theorem{ Let $L$ and $\Om$ be as in {\rm Theorem \rm 4.1} with $\Omsiga$ decomposed as
in {\rm (2.1)} and $\pa\Om_j$ given by {\rm (2.2)}. Let $\Lam\in C^\sigma(\pa\Om,\pa\D)$ be
generic, $\phi\in C^\sigma(\pa\Om,\R)$ with $0<\sigma <1$ and let
$\kappa_1,\cdots ,\kappa_{N+1}$ be the indices of $\Lam$ with respect to $\Om_1,\cdots ,\Om_{N+1}$,
respectively. Then, the Problem {\rm (5.1)} has a solution $u$ such that $u$ is continuous everywhere
on $\ov{\Om}$ except (possibly) at isolated points. The singular points of $u$ consist of a single pole of order
$\le -\kapj$ in each component $\Om_j$  where $\kapj <0$. Moreover, the value of
any such solution $u$ on an orbit $\Ga=\mathrm{arc}(p^-,p^+)$ is uniquely determined and it is given by
\[
u(\Ga)=\ep\,\frac{\Lam(p^-)\phi(p^+)-\Lam(p^+)\phi(p^-)}{i\sin(\pi\alpha)}\, ,
\]
where $\alpha\in (0,\ 1)$ is the associated jump of $\Lam$ along $\Ga$ and $\ep =1$ or $-1$.   }}

\begin{proof}
Let $Z:\, \ov{\Om^0}\,\longrightarrow\,\ov{\D}$ be the first integral of $L$ as defined in Theorem 2.3.
We can assume that $Z^{-1}(0)\cap\sigc =\emptyset$ (if not replace $Z$ by $H\circ Z$ with $H$ an appropriate
automorphism of $\D$). Hence, $Z^{-1}(0)$ consists of $N+1$ points, one point in each component $\Om_j^0$.
For $j\in\{1,\cdots ,N+1\}$, let $Z_j:\,\ov{\Om_j^0}\,\longrightarrow\,\ov{\D}$ be the
restriction of $Z$ to
$\ov{\Om_j^0}$ (given by (2.4)) and with $\pa{\Om_j^0}$ given by (2.5).
As in Theorem 2.3,  for $ l=1,\cdots ,p$ let
$z_{jl}=Z_j(C_l)\in\D$ and for $k=1,\cdots ,m$ let
$\cjk =Z(\Gamjk)\in\pa\D$. Let $\lamj=\Lam\circ Z_j^{-1}$ be defined on
$\pa\D\backslash\{ c_{j1},\ldots ,c_{jm}\}$ and let $\lam_j^0$ be the modification
of $\lamj$ as constructed in the proof of Theorem 4.1. Hence, $\lam_j^0$ is H\"{o}lder
continuous on $\pa\D$ and is defined by (4.13) or (4.14). The function
\begeq
w_j^0(z)=\frac{u\circ Z_j^{-1}(z)}{\prod_{k=1}^m(z-\cjk)^{\aljk}}\, ,\qquad z\in\D
\stopeq
satisfies the nonhomogeneous RH-Problem
\begeq
\dd{w_j^0(z)}{\ov{z}}=0\ \ z\in \D ,\quad \mathrm{Re}\left(\ov{\lam_j^0}(\zeta)u(\zeta)\right)=\rho(\zeta)\ \
 \zeta\in\pa\D\, ,
\stopeq
where
\begeq
\rho (\zeta) =\frac{\beta (\zeta)\psi(\zeta)}{\prod_{k=1}^m\abs{\zeta -\cjk}^{\aljk}}\, ,
\stopeq
with $\beta$ the alternating function defined in (4.12) and $\psi(\zeta)=\phi\circ Z_j^{-1}(\zeta)$
for $\zeta\in\pa\D$ and $\zeta\ne \cjk$.

Suppose first that $\kapj \ge 0$. Then, the general solution of Problem (5.6) is given by
\begeq
w_j^0(z)=z^{\kapj}\ei{i\gamma(z)}\left(
\mathcal{S}(\widehat{\rho})(z)+id_0+\sum_{l=1}^{\kapj}d_lz^l-\ov{d_l}\frac{1}{z^l}\right)\, ,
\stopeq
where $d_0\in\R$, $d_1,\cdots ,d_{\kapj}\in\C$ are arbitrary constants, the function $\gamma(z)$ is given
by (4.19), $\mathcal{S}$ is the Schwarz operator in $\D$ given by (4.20), and where
\begeq
\widehat{\rho}(\zeta)=\ei{\mathrm{Im}(\gamma(\zeta))}\rho(\zeta)\, .
\stopeq
We need to understand the behavior of $\mathcal{S}(\widehat{\rho})(z)$ as $z\in\D$ approaches a boundary
point $\cjk$. For this we rewrite $\mathcal{S}(\widehat{\rho})(z)$ as
\begeq
\mathcal{S}(\widehat{\rho})(z)=\frac{1}{2\pi i}\int_{\pa\D}\frac{\widehat{\rho}(\zeta)}{\zeta}\, d\zeta +
\frac{2z}{2\pi i}\int_{\pa\D}\frac{\widehat{\rho}(\zeta)/\zeta}{\zeta -z}\, d\zeta\, .
\stopeq
The function $\widehat{\rho}(\zeta)/\zeta$ can be written as
\begeq
\frac{ \widehat{\rho}(\zeta)}{\zeta}=\frac{g(\zeta)}{(\zeta -\cjk)^{\aljk}}\, ,
\stopeq
with
\begeq
g(\zeta)=\frac{\ei{\mathrm{Im}(\gamma(\zeta))}\beta(\zeta)\psi(\zeta)}{\zeta\prod_{l\ne k}
\abs{\zeta-c_{jl}}^{\alpha_{jl}}}\times\frac{(\zeta-\cjk)^\aljk}{\abs{\zeta-\cjk}^\aljk}\, .
\stopeq
With our choice of a branch of the logarithm as in the proof of Theorem 4.1, we have
\begeq\begar{ll}
\dis g(\cjk^+)=&\dis
\frac{\ei{\mathrm{Im}(\gamma(\cjk))}\beta(\cjk^+)\psi(\cjk^+)}{\cjk\prod_{l\ne k}\abs{\cjk -c_{jl}}^{\alpha_{jl}}}\,;\\
\\
\dis g(\cjk^-)=&\dis
\frac{\ei{\mathrm{Im}(\gamma(\cjk))}\beta(\cjk^-)\psi(\cjk^-)
\ei{i\pi\aljk}}{\cjk\prod_{l\ne k}\abs{\cjk -c_{jl}}^{\alpha_{jl}}}\, .
\stopar\stopeq
It follows from (5.10), (5.11), (5.12) and Lemma 5.1 that
\begeq
\mathcal{S}(\widehat{\rho})(z)=
\frac{\ei{\mathrm{Im}(\gamma(\cjk))}}{\cjk\prod_{l\ne k}\abs{\cjk -c_{jl}}^{\alpha_{jl}}}\times
\frac{M_{jk}}{(z-\cjk)^\aljk}+\Phi_0(z)\, ,
\stopeq
with
\begeq
\abs{\Phi_0(z)}\le\frac{K}{\abs{z-\cjk}^{\alpha^0_{jk}}}
\stopeq
for some $K>0$ and $0<\alpha_{jk}^0<\aljk$ and where
\begeq
M_{jk}=\frac{\ei{i\pi\aljk}\beta(\cjk^+)\psi(\cjk^+)-\beta(\cjk^-)\psi(\cjk^-)}{2i\sin(\pi\aljk)}\, .
\stopeq
Therefore, it follows from (5.8) and (5.14), that the function
\begeq
w_j(z)=w_j^0(z)\prod_{k=1}^m(z-\cjk)^\aljk
\stopeq
satisfies (for $z\in\D$ close to $\cjk$)
\begeq
w_j(z)=2z^{\kapj+1}\frac{\ei{i\gamma(z)}\ei{\mathrm{Im}(\gamma(\cjk))}}{\cjk}\times\prod_{l\ne k}
\frac{(z-c_{jl})^{\alpha_{jl}}}{\abs{\cjk-c_{jl}}^{\alpha_{jl}}}\times M_{jk}+w_j^\ast(z)\, ,
\stopeq
with $w_j^\ast$ analytic and $w_j^\ast(\cjk)=0$. Therefore,
\begeq
w_j(\cjk)=\lim_{z\in\D\,\to\, \cjk}w_j(z)=2\cjk^{\kapj}\ei{i\mathrm{Re} (\gamma(\cjk))}M_{jk}
\prod_{l\ne k}\frac{(\cjk-c_{jl})^{\alpha_{jl}}}{\abs{\cjk-c_{jl}}^{\alpha_{jl}}}\, .
\stopeq
Since $\gamma$ is given by (4.19), then
\begeq
\ei{i\mathrm{Re}(\gamma(\cjk))}=\cjk^{-\kapj}\lamj^0(\cjk)=
\cjk^{-\kapj}\beta(\cjk^+)\lamj(\cjk^+)\prod_{l\ne k}
\frac{\ov{(\cjk-c_{jl})^{\alpha_{jl}}}}{\abs{\cjk-c_{jl}}^{\alpha_{jl}}}\, .
\stopeq
Hence, using (5.20) and (5.16) in (5.19), we get
\begeq\begar{ll}
w_j(\cjk)& =2\beta(\cjk^+)\lamj(\cjk^+)M_{jk}\\
\\
& =\dis\frac{\lamj(\cjk^+)\psi(\cjk^+)\ei{i\pi\aljk}-
\beta(\cjk^+)\beta(\cjk^-)\lamj(\cjk^+)\psi(\cjk^-)}{i\sin(\pi\aljk)}\, .
\stopar\stopeq
 The expression for $w_j(\cjk)$ can be
further simplified by noticing that
\begeq
\beta(\cjk^+)\beta(\cjk^-)=(-1)^{\qjk}\quad\mathrm{and}\quad
\lamj(\cjk^-)=(-1)^{\qjk}\ei{i\pi\aljk}\lamj(\cjk^+)
\stopeq
so that
\begeq
w_j(\cjk)=(-1)^{\qjk}
\frac{\lamj(\cjk^-)\psi(\cjk^+)-\lamj(\cjk^+)\psi(\cjk^-)}{i\sin(\pi\aljk)}\, .
\stopeq
In terms of the values of original functions $\Lam$ and $\phi$ at the ends $\pjkm$ and $\pjkp$
of the orbit $\Gamjk$, we have
\begeq
w_j(\cjk)=(-1)^{\qjk}
\frac{\Lam(\pjkm)\phi(\pjkp)-\Lam(\pjkp)\phi(\pjkm)}{i\sin(\pi\aljk)}\, .
\stopeq

In summary, going back to the original problem, we see that
in the case $\kapj \ge 0$, the function $u_j$ defined in $\ov{\Om_j}$ by
\begeq
u_j(x,y)=\left\{\begar{ll}
w_j(Z_j(x,y)) &\ \ (x,y)\in\ov{\Om_j^0}\, ,\\
w_j(z_{jl}) &\ \ (x,y)\in\intc{l},\ \ l=1,\cdots ,p\, ,
\stopar\right.\stopeq
is continuous on $\ov{\Om_j}$,  satisfies
\begeq
Lu=0\ \ \mathrm{in}\ \Om_j\, ,\quad \mathrm{Re}(\ov{\Lam}(t)u(t))=\phi(t)\ \
t\in\pa{\Om_j}\cap\pa\Om\, ,
\stopeq
and is constant on each arc $\Gamjk$ along which it is uniquely determined by
$u(\Gamjk)=w_j(\cjk)$ given by (5.24).

Now we consider the case $\kapj <0$. The solution
\begeq
w_j^0(z)=z^{\kapj}\ei{i\gamma(z)}\mathcal{S}(\widehat{\rho})(z)
\stopeq
of Problem (5.6) is meromorphic with a single pole at $z=0$ of order
$\le -\kapj$. The function defined in $\ov{\Om_j}\backslash Z_j^{-1}(0)$ by
\begeq
u_j(x,y)=\left\{\begar{ll}
w_j(Z_j(x,y)) &\ \ (x,y)\in\ov{\Om_j^0}\backslash Z_j^{-1}(0)\, ,\\
w_j(z_{jl}) &\ \ (x,y)\in\intc{l},\ \ l=1,\cdots ,p\, ,
\stopar\right.\stopeq
is continuous with a pole of order $\le -\kapj$ at the point $Z_j^{-1}(0)$
and solves problem (5.26) in $\Om_j\backslash Z_j^{-1}(0)$. Furthermore, its
value on each orbit $\Gamjk$ is again given by (5.24).

These various functions $u_j$ defined in $\Om_j$ yield a well defined continuous function $u$
on the whole open set $\Om$ (except at isolated poles). Furthermore,  it is easily verifiable that $u$ solves
the RH-Problem (5.1).
\end{proof}

{\remark{ The RH-Problem (5.1) could admit continuous solutions throughout $\Om$ even when
some of the indices $\kapj <0$, provided that the boundary functions $\Lam$ and $\phi$
satisfy ($-\kapj -1$)-conditions on $\pa\Om_j$. These conditions simply mean that $\mathcal{S}(\widehat{\rho})$
has a zero of order $\ge -\kapj$ at $z=0$ (where $\widehat{\rho}$ is given by (5.9)). That is,
$\dis\int_{\pa\D}\widehat{\rho}(\zeta)\zeta^{-s}d\zeta =0$ for $s=1,\cdots ,-\kapj$. In terms of integrals
over $\pa\Om_j$, we get the conditions
\begeq
\sum_{k=1}^m\int_{\Ajk}\!\!\widehat{\rho}(Z_j(t))\frac{dZ_j(t)}{Z_j^s(t)}=0\quad s=1,\cdots ,-\kapj\, .
\stopeq }}

\medskip

The last problem we consider in this paper is
\begeq\left\{\begar{ll}
Lu=f & \ \ \mathrm{in}\ \ \Om\, ,\\
\mathrm{Re}\left(\ov{\Lam}\, u\right) =\phi &\ \ \mathrm{on}\ \ \pa\Om\, .
\stopar\right.\stopeq
We assume that $\Lam$ and $\phi$ are as in Theorem 5.2 and $f\in C^\infty(\ov{\Om})$ has
zero periods on the closed orbits of $L$. By this, we mean that if we consider the
dual differential form $ \omega =Bdx-Ady$ of the vector field $L$, then there exists a $C^\infty$-differential form
$\eta$ in $\ov{\Om}$ such that $fdx\wedge dy=\omega\wedge\eta$, and we say that $f$ has {\it{zero periods}} on
the closed orbits of $L$ if
\begeq
\int_C\eta =0\quad\mathrm{for\ every\ closed\ orbit}\ \ C\subset\sigc\, .
\stopeq
It should be noted that this condition depends only on the function $f$ and not on the particular
choice of its representative $\eta$ (any other representative is of the form $\eta+g\omega$ for some
function $g$). It is proved in {\cite{Mez-JDE}} that condition (5.31) is a necessary and sufficient
condition for the global solvability  of $L$ on $\Om$.

Now we can state our result regarding Problem (5.30).

{\theorem{ Let $L\, ,\  \Om\, ,\ \Lam\, ,\ \phi$ be as in {\rm Theorem 5.2}. Let $\kapj$
be the index of $\Lam$ relative to the component $\Om_j$ of $\Omsiga$ and $f\in C^\infty(\ov{\Om})$
satisfies condition {\rm (5.31)}. Let $\kappa_{j_1},\,\cdots\, ,\kappa_{j_n}$ be the negative
indices of $\Lam$. Then, for every $l=1,\,\cdots\, ,n$ there exists a point
$s_{j_l}\in\Om_{j_l}\backslash\siginf$ such that Problem {\rm (5.30)} has a solution $u$ that is
H\"{o}lder continuous on $\ov{\Om}\backslash\{ s_{j_1},\,\cdots\, ,\, s_{j_n}\}$
and $C^\infty$ in $\ov{\Om}\backslash\siginf\cup\{ s_{j_1},\,\cdots\, ,\, s_{j_n}\}$.
Furthermore, for each $s_{j_l}$ there exist an open set $\mathcal{O}_{j_l}$, with
$s_{j_l}\in\mathcal{O}_{j_l}\subset\Om$, and a function $m_{j_l}\in C^\infty(\mathcal{O}_{j_l})$ such
that $u-m_{j_l}$ has a pole of order $\le -\kappa_{j_l}$ at $s_{j_l}$.
}}

\begin{proof}
We know from {\cite{Mez-JDE}} that for a given $f\in C^\infty(\ov{\Om})$ satisfying (5.31), there
exists a function $v$, H\"{o}lder continuous on $\ov{\Om}$ and in $C^\infty(\ov{\Om}\backslash\siga)$
such that $Lv=f$. Theorem 5.2 implies that the RH-Problem
\[
Lw=0\ \ \mathrm{in}\ \Om\, ,\qquad \mathrm{Re}(\ov{\Lam}\, w)=\phi -\mathrm{Re}(\ov{\Lam}\, v)
\]
has a solution $w$ such that $w$ is
H\"{o}lder continuous on $\ov{\Om}\backslash\{ s_{j_1},\,\ldots\, ,\, s_{j_n}\}$,
is $C^\omega$ in $\ov{\Om}\backslash\siginf\cup\{ s_{j_1},\,\ldots\, ,\, s_{j_n}\}$,
 and  has a pole of order $\le -\kappa_{j_l}$
at the points $s_{j_l}$. The function $u=w+v$ is then the sought solution of Problem (5.30).
\end{proof}

 \bibliographystyle{amsplain}

 \end{document}